\documentclass[11pt]{amsart}
\usepackage{latexsym}

\usepackage{amsmath,amssymb,amsfonts,array,graphicx}

\def\CC{\mathbb{C}}

\def\QQ{\mathbb{Q}}
\def\RR{\mathbb{R}}
\def\ZZ{\mathbb{Z}}
\def\SS{\mathbb{S}}
\def\PP{\mathbb{P}}
\def\delbar{\overline{\partial}}

\def\cal{\mathcal}

\author{Klaus Mohnke}
\title{How to (Symplectically) Thread the Eye of a (Lagrangian) Needle}
\address{Klaus Mohnke, Institut f\"ur Mathematik, Humboldt--Universit\"at zu
Berlin}
\thanks{Supported by
Deutsche Forschungsgemeinschaft DFG (Mo 843/1-2) and the American
Institute of Mathematics, Palo Alto}
\email{mohnke@mathematik.hu-berlin.de}
\date{}
\newtheorem{theorem}{Theorem}
\newtheorem{definition}[theorem]{Definition}
\newtheorem{proposition}[theorem]{Proposition}
\newtheorem{lemma}[theorem]{Lemma}
\newtheorem{remark}[theorem]{Remark}

\dedicatory{To Janett}
\begin{document}
\bibliographystyle{plain}

\begin{abstract}
We show that there exists no Lagrangian embeddings of the
Klein bottle into $\CC\PP^{2}$. Using the same techniques we also give a new proof
that any Lagrangian torus in $\CC\PP^2$ is smoothly isotopic to the Clifford torus.
\end{abstract}

\maketitle

\section{Lagrangian Embeddings in $\CC^{2}$}

The topology of closed Lagrangian embeddings into $\CC^{n}$ (see
\cite{Arnold:firststeps}) is still an elusive problem in
symplectic topology. Before Gromov invented the techniques of
pseudo--holomorphic curves it was almost intractable and the only
known obstructions came from the fact that such a submanifold has
to be totally real. Then in \cite{Gromov} he showed that for any
such closed, compact, embedded Lagrangian there exists a
holomorphic disk with boundary on it. Hence the integral of a
primitive over the boundary is different from zero and the first
Betti number of the Lagrangian submanifold cannot vanish,
excluding the possibility that a three--sphere can be embedded
into $\CC^{3}$ as a Lagrangian. A further analysis of these
techniques led to more obstructions for the topology of such
embeddings in \cite{Polterovich} and \cite{Viterbo:lagtori}.

For $\CC^{2}$ the classical obstructions restrict the classes
of possible closed, compact surfaces which admit Lagrangian embeddings
into $\CC^{2}$ to the torus and connected sums
of the Klein bottle with oriented surfaces of even genus.
There are obvious Lagrangian embeddings of the torus
(e.g.~the Clifford torus $\SS^{1}\times
\SS^{1}\subset\CC\times \CC$) and not so obvious ones for the connected sums
except for the Klein bottle (see
\cite{Audin/Lalonde/Polterovich}). One may further
ask which topological types of embeddings may be realized as
a Lagrangian embedding. There has been a partial answer to that in
\cite{Luttinger} and an announcement of a proof by Hofer and Luttinger
of the following
\begin{theorem}[Hofer/Luttinger]\label{unknot}
Any Lagrangian embedding of the torus into $\CC\PP^2$ is smoothly isotopic
to the Clifford torus $\SS^1\times \SS^1\subset\CC\times\CC$.
\end{theorem}

Here we show that the same circle of ideas, namely the study of holomorphic
curves with respect to a certain singular almost complex structure,
solves the question for the Klein bottle:

\begin{theorem}\label{kleinbottle}
There is no Lagrangian embedding of the Klein bottle into $\CC\PP^{2}$.
\end{theorem}

Both questions were contributed to Kirby's problem list
\cite{Kirby:problems} by Eliashberg. There have been several
attempts to attack this problem (\cite{Audin/Lalonde/Polterovich,
Lalonde, Mohnke:Lagrange} et.al.). Independently, Nemirovski
proved in \cite{Nemirovski:kleinbottle} a more general version of
Theorem~\ref{kleinbottle} using completely different methods from
complex analysis.

The constructions described in this paper will yield
the following result:

\begin{theorem}\label{threadings}
Let $L\subset \CC\PP^2$ be a closed Lagrangian embedding of a flat
surface, i.e.~either the Klein bottle or the $2$--torus. Then
there is an almost complex structure $J$ compatible to $\omega$
with the following property: There are  two smooth mappings
$D,E:\Delta\to \CC\PP^2$ of the unit disk with
$D(\partial\Delta),E(\partial\Delta)\subset L$ and three embedded
$J$--holomorphic spheres $F,G,H:\CC\PP^1\to \CC \PP^2\setminus L$
whose fundamental classes represent the same  generator of
$H_2(\CC\PP^2;\ZZ)\cong\ZZ$ with the following algebraic
intersection numbers:
\begin{align}\label{pairintersect}
F\cdot D &= 1, \quad F\cdot E = 0\nonumber\\
G\cdot D &= 0, \quad G\cdot E = 1\nonumber\\
H\cdot D &= 0, \quad H\cdot E = 0.\nonumber
\end{align}
\end{theorem}

Let us first see how the main result follows and how we obtain a
new proof of the unknottedness of Lagrangian tori.

\begin{proof}[Proof of Theorem~\ref{kleinbottle}]
We claim that the span of $[\partial D],[\partial E]\in
H_1(L;\QQ)$ is two--dimensional. Assume for some pair of integers
$(k,l)\neq (0,0)$
\begin{equation}
k[\partial D] = l[\partial E]
\end{equation}
as {\em integer} homology classes. Hence we can find a $2$--chain
$C$ in $L$ such that $kD+ C + (-lE)$ forms a $2$--cycle in
$\CC\PP^2$. We test it against $F,G$ and $H$ which yields the same
number since these represent the same homology class. By the
assumption, $F\cdot C = G\cdot C= H\cdot C = 0$ and therefore
\begin{equation}
   \begin{aligned}
         (kD+ C + (-lE))\cdot F &= k\\
         (kD+ C + (-lE))\cdot G &= -l\\
         (kD+ C + (-lE))\cdot H &= 0
   \end{aligned}
\end{equation}
Thus $k=l=0$, verifying the statement. Hence, $b_1(L)\ge 2$, which
excludes that $L$ was the Klein bottle.
\end{proof}

\begin{proof}[Proof of Theorem~\ref{unknot}]
By a result of Gromov in \cite{Gromov} there exists a symplectic
isotopy of $(\CC \PP^2, \omega)$ which maps the $J$-holomorphic
spheres $F,G,H$ to projective lines. Hence the image of $L$ under
this map, which we denote by abuse of notation by $L$ again, is a
Lagrangian in $\CC^2$ in the complement of $2$ complex lines,
without loss of generality, $0\times\CC, \CC\times 0\subset
\CC^2$. The lines correspond to $F$ and $G$, respectively. Hence
the intersection numbers with the symplectic disks show that the
induced homomorphism
$$
H_1(L)\longrightarrow H_1(\CC^2\setminus(0\times\CC\cup \CC\times 0))
$$
is an isomorphism. But
$\CC^2\setminus(0\times\CC\cup \CC\times 0)\cong T^\ast T^2$
symplectomorphically.
The standard Clifford torus $\SS^1\times \SS^1\subset \CC^2$ is thereby
mapped to the zero section of the cotangent bundle
$O_{T^2}\subset T^\ast T^2$. Moreover, the Lagrangian embedding
$L\subset T^\ast T^2$
induces an isomorphism $H_1(L)\cong H_1(T^2)$. By a result of Eliashberg and
Polterovich \cite{Eliashberg/Polterovich:unknottedness} it follows that
$L\subset T^\ast T^2$ is smoothly isotopic to $O_{T^2}$ in $T^\ast T^2$
and hence the original $L\subset\CC^2$ is smoothly isotopic
to the Clifford torus.
\end{proof}

\begin{remark}
\begin{enumerate}
\item The original idea of the proof goes back to earlier work
\cite{Mohnke:Lagrange}, where we observed that there is no
Lagrangian embedding of the Klein bottle into the complement of a
complex line in $\CC^2$ nontrivially linking it. The proof
presented here would show that an embedding into $\CC^2$ is
symplectically isotopic to one in the complement of a complex
line.  The nontrivial linking is  measured by the disk $D$. Thus
the existence of a Lagrangian embedding into $\CC\PP^2$ would
contradict the previous result. This idea gave the title for this
exposition.
\item Notice that such a result is bound to dimension $4$ not only for
the lack of homological intersection of pairs of holomorphic
curves: In dimension $6$ and higher there {\em are} Lagrangian
embeddings of manifolds with first Betti number $1$ (see
\cite{Mohnke:Lagrange}). Furthermore, the proof of the result does
not seem to apply directly to the cases of Lagrangian embeddings
of (non--orientable) surfaces of higher genus. Thus one could
wonder if it is true at all in these cases.
\item The argument in \cite{Eliashberg/Polterovich:unknottedness} fails
to provide an isotopy through Lagrangian embeddings.
Hence this stronger question
remains open in both situations, $L\subset \CC^2$ as well as
$L\subset T^\ast T^2$.
\end{enumerate}
\end{remark}

{\em Acknowledgements.} I am indebted to Fr\'{e}d\'{e}ric
Bourgeois, Kai Cieliebak, Yakov Eliashberg, John Etnyre and Helmut
Hofer for fruitful discussions and their generous help  to make my
idea for the proof work. I want to express especially warm
appreciation for the encouraging support I received from Helmut
Hofer, Abbas Bahri and Sagun Chanillo in a difficult moment. I was
deeply moved and without that, this work would most likely have
not been carried out. I like to thank the Department of
Mathematics at Stanford University for its hospitality during my
visit when most of the result was accomplished. Moreover, I would
like to thank the SFB 288 and Humboldt--Universit\"at zu Berlin,
Germany, the Lorentz Center at Leiden Universiteit, The
Netherlands, and the American Institute of Mathematics for their
hospitality and support. I owe special thanks to Helga Baum, Paul
Biran, Slava Charlarmov, Tobias Ekholm, Siddharta Gadgil,
Hansj\"org Geiges, Wilhelm Klingenberg, Fran\c{c}oise Lalonde,
Dusa McDuff and Dietmar Salamon.

\section{Geometric set--up and data}\label{geometry}

Let us first introduce and discuss the geometric data which
will be used to define holomorphic curves.

\subsection{Topology and Dynamics of the Unit Cotangent
Bundles}\label{dynamics}

Throughout let $L$ be either the Klein bottle $K^{2}$ or the
$2$--torus $T^{2}$ with a flat Riemannian metric $g$ fixed. Denote
by $\theta\in\Omega^1(T^\ast L)$ the canonical, or Liouville form.
Remember that its differential $d\theta$ is a symplectic form on
$T^\ast L$. The restriction to the co-sphere bundle $\lambda:=
\theta|_{U^\ast L}$, $U^\ast L:= \{\alpha\in T^\ast L\mid
|\alpha|_g = 1\}$ is a contact form We denote by
$\xi:=\ker\lambda$ the corresponding contact distribution. The
flow of its {\em Reeb vector field} $R_\lambda$, given by
\begin{equation}
   \begin{aligned}
       d\alpha(R_\lambda) &= 0\\
       \alpha(R_\lambda) &=  1,
   \end{aligned}
\end{equation}
coincides with the geodesic flow of the chosen metric. Hence,
closed orbits of the Reeb flow, in short {\em Reeb orbits},
correspond exactly to closed {\em oriented} geodesic Notice that
the Reeb flow preserves $\lambda$, and therefore induces
symplectic isomorphisms of the symplectic vector spaces $(\xi,
d\lambda)$. We will not distinguish (oriented) geodesics and their
corresponding closed Reeb orbits by different notation. From the
context it will be clear wether we denote by $\gamma$ either of
them. It will not necessarily be simple. By $-\gamma$ we denote
the geodesic with the opposite orientation or its corresponding
Reeb orbit. Notice that the latter is geometrically different
(actually disjoint) from $\gamma$.

The complement of the zero section $0_L$ in $T^\ast L$ is symplectomorphic to
the symplectization of $U^\ast L$:
\begin{equation}\label{ident}
 \begin{aligned}
  (T^\ast L\setminus 0_L, d\theta) &\cong (\RR\times U^\ast L
  d(e^r\lambda)),\\
  v &\mapsto (\frac{v}{|v|}, \log(|v|)),
 \end{aligned}
\end{equation}
where $r\in \RR$ denotes the parameter of the second factor.
An almost complex structure $J$ on  $(\RR\times U^\ast L, d(e^r\lambda))$
is called {\em compatible} (to $\lambda$) if
\begin{equation}
   \begin{aligned}
         J(\frac{\partial}{\partial t}) &= R_\alpha\\
         J(\xi) &= \xi\\
         d\lambda(., J.) &> 0 \mbox{ and symmetric }\\
         J &\quad  \mbox{invariant under translation}.
   \end{aligned}
\end{equation}
Therefore, $J$ is compatible to the symplectic structure:
$d(e^r\alpha)(.,J.)$ is a Riemannian structure (a symmetric and
positive definite bilinear form). It is a cone over $U^\ast L$.
Most importantly, cylinders over  Reeb orbits are
$J$--holomorphic. Finally, we call a complex structure $J$ on
$T^\ast L$ compatible (to $d\theta$ and $\lambda$) if
$d\theta(.,J.)$ defines a Riemannian metric and if $J$ coincides
with a $\lambda$--compatible almost complex structure on the
complement $T^\ast L \setminus D_\rho^\ast L$ of $D_\rho^\ast
L:=\{\alpha\in T^\ast L\mid |\alpha|_g\ge \rho\}$ under the above
identification (\ref{ident}), $T^\ast L \setminus D_\rho^\ast
L\cong (\log\rho,\infty)\times U^\ast L$.

On the other hand the metric $g$ on $L$ defines an almost
complex structure $J_g$ on $T^\ast L$, which is compatible to
$d\theta$ in the following way:
we define a complex structure on $T^\ast L\oplus TL$ via
$$
J_g(\alpha,v) := (-g(v,.), g^{-1}\alpha).
$$
Then we use the Levi-Civita connection to globally split $T(T^\ast
L)\cong T^\ast L \oplus TL$ into horizontal and vertical part with
respect to the fibration $T^\ast L\longrightarrow L$. $J_g$ is not
translational invariant, and hence not compatible to $\lambda$.
Instead, we may use it to define a $\lambda$-compatible almost
complex structure $J$: We set $J(x):= J_g(x)$ for $x\in U^\ast L$
and extend it to $\RR\times U^\ast L$ via translation.

Given a geodesic $\gamma=\gamma(t)$ in $L$, the map
$f_\gamma:\CC\to T^\ast L$ given by
$$
f_\gamma(s+\mbox{i}t):= sg(\dot{\gamma}(t),.)
$$
is $J_g$--holomorphic. Since the image in the complement of $0_L$ consists
of two cylinders over the closed Reeb orbits, $\RR\times \pm\gamma$,
these are also $J$--complex curves.
But they will become holomorphic only for a different
parameterization.

Now we use a cut--off function $\varphi=\varphi(|\alpha|)$ to
interpolate between $J_g$ and $J$ and obtain a
$d\theta$-compatible almost complex structure $J^-$ on $T^\ast L$
which coincides with $J_g$ in a neighborhood of the zero section
$0_L$ and with $J$ on $[-1,\infty)\times T^\ast_1 L$ with respect
to  (\ref{ident}). Under these assumptions, the image of
$f_\gamma$ is still a complex curve. By abuse of notation we will
denote by $f_\gamma$ its $J^-$--holomorphic parameterization.

Closed geodesics on $L$ typically arise in $1$--dimensional
families. However the corresponding Reeb orbits are
non--degenerate in the sense of Morse--Bott: any eigenvector of
eigenvalue $1$ of the linearized Poincar\'{e} return map
restricted to $\xi$ (which is transversal to it) arises as the
deformation of the closed Reeb orbit through a family of such (see
\cite{Hofer/Wysocki/Zehnder:Degeneracies}). Hence we may use the
results of \cite{Hofer/Wysocki/Zehnder:Degeneracies}, or
\cite{Bourgeois:Morse-Bott, Bourgeois/Eliashberg:Compactness} in
our argument. The advantage over perturbing $\lambda$ such that
there are only non--degenerate closed Reeb orbits in the strong
sense is a clearer structure of the proof.

If $L$ is the flat Klein bottle there are also  isolated
geodesics, two simple ones and their odd multiples which are
parallel (but not homotopic). They  are all non-degenerate in the
stronger sense. But notice that their even multiple covers are
elements of $1$--dimensional families of geodesics as before. Thus
let us call the isolated geodesics {\em odd} and all those which
arise in families {\em even}.

Now there are three indices one can assign to a geodesic $\gamma$,
which is non-degenerate in the sense of Morse--Bott. First of all
its {\em Morse index}, $\mbox{index}(\gamma)$: it is the maximal
dimension of a linear subspace in the space of vector fields at
$L$ along $\gamma$ on which the Hessian of the length functional
on loops in $L$ is negative definite.

For the other two indices we  fix a trivialization of the
symplectic vector bundle $(T(T^\ast L),d\theta)$. With that choice
any loop  $\gamma:\SS^1\to \subset L$ defines a loop of Legendrian
subspaces $\Gamma(x) = T_{\gamma(x)}L\subset \CC^2$ and we may
assign to the geodesic $\gamma$ the {\em Maslov index}
$\mu(\gamma)$. It depends only on the homology class $[\gamma]$ of
a loop in $L$ and hence defines a homomorphism
$$
\mu:H_1(L;\ZZ)\longrightarrow \ZZ.
$$
Recall that $\mu(\gamma)$ is even if and only if $\gamma^\ast TL$
is orientable.
 Notice that for $L=T^2$ we may choose the
trivialization induced by a trivialization of $TL\to L$. In this
case all Maslov--indices vanish.

On the other hand this trivialization induces one of the
symplectic bundle $(\xi, d\lambda)$ via
$$
\xi\oplus \RR R_\lambda\oplus \RR \frac{\partial}{\partial r}
\cong T(T^\ast L)
$$
on $U^\ast L$. Since the linearization of the Reeb flow induces
symplectic isomorphisms of $(\xi, d\lambda)$ any closed Reeb orbit
defines a path of symplectic $2\times 2$--matrices, starting at
the identity. Thus we can assign the Conley--Zehnder index
$CZ(\gamma)$  to the closed Reeb orbit $\gamma$, which, in our
situation, could be a half--integer (see
\cite{Robbin/Salamon:maslovforpaths}).

Viterbo found  a striking relation between them:
\begin{lemma}[\cite{Viterbo:lagtori}]\label{viterbo}
\begin{equation}
CZ(\gamma) = \mu(\gamma) + \mbox{index}(\gamma) +
\frac{1}{2}\dim(\gamma).
\end{equation}
Here $\dim(\gamma)$ is the dimension of the family of geodesics in
which $\gamma$ arises.
\end{lemma}
Notice that Viterbo gave this formula in the non-degenerate case
in the strong sense, but this identity follows from his result and
the result of \cite{Robbin/Salamon:maslovforpaths}.

\subsection{Almost complex structures with cylindrical ends}\label{concaveend}

Let $L\subset \CC\PP^2$ be a Lagrangian embedding. There is a
(pseudo-convex) neighborhood $U$ of $L\subset \CC\PP^{2}$ which is
symplectomorphic to the unit disk bundle $D_1^\ast L\subset T^\ast
L$ of $L$ after rescaling the flat metric on $L$. Fix an almost
complex structure $J_0$ on $\CC\PP^2$ which is compatible with the
symplectic structure $\omega$, coincides with $J^-$ on $U$ using
identification (\ref{ident}). These data define an
$\omega$--compatible almost complex structure $J^+$ on
$\CC\PP^2\setminus L$  which coincides with $J_0$ on
$\CC\PP^2\setminus U$ and with $J$ on $U\setminus L$ again using
(\ref{ident}). $\CC\PP^2\setminus L$ is said to have a concave
end, and $J^+$ is compatible to $\omega$ and $\lambda$. Notice
that the Riemannian metric defined by $\omega(.,J^+.)$ has a
singularity which is a cone over the unit cotangent bundle $U^\ast
L$ equipped with the induced metric.

Now we deform $J_0$ into a family $J_\tau$ of $\omega$--compatible
almost complex structures to ''approximate'' $J^+$. For a
convenient description we define spaces $X_\tau$ via
$$
X_\tau := U \cup_{\partial U\cong \{-\tau\}\times U^\ast L}
[-\tau,0]\times U^\ast L \cup_{\{0\}\times U^\ast L\cong
\partial(\CC\PP^2\setminus U)}\CC\PP^2\setminus U.
$$
Then $J_\tau$ coincides with $J^+$ on $\CC\PP^2\setminus U$, with
$J^-$ on $U$ and with $J$ on the cylinder. It is compatible to a
symplectic structure $\omega_\tau$ which is given by $\omega$ on
$\CC\PP^2\setminus U$, by $d(e^r\alpha)$ on $[-\tau,0]\times
U^\ast L$,and by the rescaled structure $e^{-\tau}d\theta$ on $U$.
Notice that there is a diffeomorphism $\Phi_\tau:\CC\PP^2\cong
X_\tau$ simply given by squeezing the long neck to its original
size. With respect to this identification the cohomology classes
of $[\omega_\tau]$ are constant
 $[\Phi_\tau^\ast\omega_\tau] = [\omega]\in H^2(X_\tau)\cong H^2(\CC\PP^2)$.
 By Moser's argument they
are all symplectomorphic by  diffeomorphisms which are supported in $U$.
Hence $J_\tau$ can all be considered as almost complex structures on $\CC\PP^2$
tamed by $\omega$. We remark that this is not crucially important for the
result presented here but added for the convenience of the reader.

\section{Holomorphic Curves}

In this  section we collect the parts of the theory necessary for
the proof of our main result. This is a  version of what is known
about punctured holomorphic curves in symplectic manifolds with
concave and convex ends (see \cite{EGH:SFT} Chapters~1.1.--1.6,
\cite{BEHWZ,Bourgeois:Morse-Bott, Cieliebak/Mohnke:punc}) reduced
to what is actually needed in the argument. The geometric set--up,
the compactness result and the index formula are specialized to
the situation at hand.

\subsection{Punctured holomorphic spheres}
Families of compact $J_\tau$--holomorphic curves in symplectic
manifolds with ends split into  $J^\pm$ and $J$--holomorphic
curves as $\tau\to\infty$ (see below). As part of the limit appear
so--called punctured pseudo-holomorphic curves. We may restrict
our considerations to spheres:
\begin{definition}
A {\em punctured holomorphic sphere} in $X=\CC\PP^2\setminus L$,
$T^\ast L$ or $\RR\times U^\ast L$ is a punctured sphere
$(\CC\PP^1, \overline{z}_1,\dots, \overline{z}_{\overline{p}},
\underline{z}_1,\dots, \underline{z}_{\underline{p}})$ together
with a map $f$ holomorphic with respect to $J^+$, $J^-$ or $J$
$$
f:\dot{\Sigma}:= \CC\PP^1\setminus\{\overline{z}_1,\dots,
\overline{z}_{\overline{p}}, \underline{z}_1,\dots,
\underline{z}_{\underline{p}}\}\longrightarrow X.
$$
For each of the punctures $\overline{z}_i$ ($\underline{z}_j$)
there exists a closed Reeb orbit $\gamma(t)=\overline{\gamma}_i$,
($\underline{\gamma}_j$), $\dot{\gamma}(t) = R(\gamma(t))$, with
the following property: In local complex coordinates $(s,t)\in
\RR_+\times \SS^1\mapsto z_i + e^{-(s+\mbox{i}t)}$ ($(s,t)\in
\RR_-\times \SS^1\mapsto z_i + e^{s+\mbox{i}t}$)around $z_i$, the
following limit is uniform with respect to $t$
\begin{equation}\label{asymptotics}
   \lim_{s\to\pm\infty}f(s,t) = \gamma(\frac{T}{2\pi}t),
\end{equation}
where $T$ is the period (or action) of $\gamma$
(see
\cite{Hofer/Wysocki/Zehnder:Degeneracies,Hofer/Wysocki/Zehnder:Asymptotics}).
\end{definition}
In the aforementioned papers the authors show that punctured
$J$--holomorphic curves are asymptotic to cylinders over the
corresponding  closed Reeb orbits near their punctures with
exponential decays with respect to $s$.

The asymptotic behavior allows us to compactify $u$ to a smooth
map
$$
\overline{f}: \overline{\Sigma}\longrightarrow
\overline{\CC\PP^2\setminus L}
$$
or $\overline{T^\ast L}$, $\overline{\RR\times U^\ast L}$,
respectively, where $\overline{\Sigma}$ is the compactification of
$\dot{\Sigma}$ by  circles (one for each puncture) and
$\overline{\CC\PP^2\setminus L}$ ($\overline{T^\ast L}$,
$\overline{\RR\times U^\ast L}$) are the compactifications  by
copies of  $U^\ast L$ each time using the cylindrical structures
near the ends. The topology of these compactifications and the
asymptotics mentioned above are such that the map
$\pi:\overline{\CC\PP^2\setminus L}\to \CC\PP^2$ induced by the
projection $\pi:_1^\ast L\to L$ is smooth. Hence, finally, we
obtain a smooth map
\begin{equation}\label{compactify}
\begin{aligned}
\hat{f} &: \overline{\Sigma}\longrightarrow \CC\PP^2\\
\hat{f} &= \pi\circ \overline{f}
\end{aligned}
\end{equation}
with $\hat{f}(\partial \overline{\Sigma})\subset L$.

Notice that the $J^-$--holomorphic curves $f_\gamma$ in $T^\ast L$
associated to a geodesic $\gamma$ are examples of two--punctured
holomorphic curves in the cotangent bundle. The next lemma
summarizes all facts we will need here.

\begin{lemma}\label{facts}
\begin{enumerate}
\item There are no
      one--punctured $J^-$--holomorphic spheres in $T^\ast L$, and no
      one--punctured $J$--holomorphic spheres in $\RR\times U^\ast L$.
\item Any holomorphic sphere with two punctures in $T^\ast L$ which
      intersects the zero section, is of the form $f_\gamma$ for a
      closed geodesic $\gamma$ as introduced above.
\item The homological intersection index with $\ZZ_2$--coefficients between
      $f_\gamma$ and $0_L$ is $1$ if  $\gamma$ is an odd closed geodesic.
\item For each odd geodesic $\gamma$ there is exactly one two--punctured
      $J^-$--holomorphic sphere in $T^\ast L$ which is asymptotic to
      $+\gamma$ ($-\gamma$) at one of its punctures, namely $f_\gamma$.
\end{enumerate}
\end{lemma}
\begin{remark}
The picture is somewhat incomplete, though it suffices for our
purposes. We expect that for each rational slope the
two--punctured spheres with asymptotics parallel to that slope
foliate $T^\ast L$  similar to the $J_g$--holomorphic shifts of
$f_\gamma$ by harmonic $1$--forms. However, the shifts are not
$J^-$--complex anymore which prevents this from being an
elementary fact.
\end{remark}
\begin{proof}
The first statement is obvious since there are no contractible
closed geodesics for a flat structure on $L$.

Let $f$ be a two--punctured holomorphic curve  with asymptotics
$\gamma',\gamma''$ which intersects $0_L$. Since geodesics
parallel to $\gamma'$ and $\gamma''$ foliate $L$ there is a
$f_\gamma$ with $\gamma$ parallel to them which intersects $f$
in $f\cap 0_L$.
Unless $f$ coincides with $f_\gamma$ (in which case we are done)
they intersect in a discrete set of points with a finite positive
(algebraic) intersection index in each of them
(see \cite{McDuff:intersections, Micallef/White}. Pick one of these
points then it will persist if we change $\gamma$ slightly. Thus
we may assume that $\pm\gamma \neq \gamma',\gamma''$.

There is an involution on $L$ which preserves all geodesics
parallel to $\gamma$ but reverses their orientations. That gives
rise to a canonical involution on $T^\ast L$ which preserves
$J^-$. Thus we obtain a curve $\tau(f)$ with asymptotics
$-\gamma', -\gamma''$. Since  $\pm\gamma \neq \gamma',\gamma''$ we
obtain a $2$--cycle in $\overline{T^\ast L}$ via
$\hat{\tau(f)}+\hat{f}-\hat{f_{\gamma'}}-\hat{f_{\gamma''}}$ which
non-trivially intersects $f_\gamma$. This is impossible, since the
latter can be moved away from the zero section through shifting it
by a nowhere vanishing $1$--form and $H_2(T^\ast L;\ZZ_2)\cong
H_2(L;\ZZ_2)$ is, of course, generated by $[0_L]$.

For the third statement we perform an explicit deformation of
$0_L$ which intersects $f_\gamma$ once transversally. There is a
one form $\beta\in\Omega^1(L)$ which vanishes on the geodesics
parallel to $\gamma$ and intersects with $0_L$ in one geodesic
$\gamma'$ perpendicular to $\gamma$. Thus $\mbox{image}(\beta)$
and $f_\gamma$ intersect in $\gamma\cap\gamma'$, and this
intersection can be made easily transversal.

To prove the last assertion notice first of all that the two
simple odd closed geodesics of the flat Klein bottle are not
homotopic. Hence, if a two--punctured holomorphic curve $f$ is
asymptotic to an odd geodesic $\gamma$ at one puncture it has to
be asymptotic to $-\gamma$ at the other. Therefore,
$\hat{f}-\hat{f_\gamma}$ forms a $2$--cycle in $\overline{T^\ast
L}$ which intersects $0_L$ homological trivially. Thus, due to the
previous fact the intersection index of $f$ and $0_L$ with
$\ZZ_2$--coefficients is non-trivial. From the second assertion we
conclude that $f=f_\gamma$.
\end{proof}

\subsection{The moduli of punctured $J$--holomorphic spheres}

$J$--holomorphic curves described in the last section typically
arise in families neglecting reparameterizations. Let us study
deformations of such a $J$--holomorphic curve $f$. Notice that
the Reeb orbits to which the curves are asymptotic to at
their punctures may also vary within the family of Reeb orbits in which they
arise. Since all closed Reeb orbits (or geodesics) are non-degenerate in the 
sense of Morse--Bott the curves are asymptotic to cylinders over closed Reeb orbits near their
punctures with exponential decays with respect to $s$ (see
\cite{Hofer/Wysocki/Zehnder:Degeneracies,
Hofer/Wysocki/Zehnder:Asymptotics}). Hence tangencies to deformations
are elements in the kernel of a Fredholm operator $\delbar_f$, which is defined
using the linearization of the Cauchy--Riemann equation acting on sections
of the pull-back, $f^\ast T\CC\PP^2$ with similar asymtptotics
additionally taking into account the variation of the Reeb orbits
(see \cite{Bourgeois:Morse-Bott}, Paragraph~2.4
for details).

If $\delbar_f$ is
surjective the virtual dimension actually coincides with the
dimension of all possible deformations modulo
reparameterizations. Such a punctured holomorphic curve is called
regular. In this case the dimension is equal to the index of
$\delbar_f$ which is given by the following formula (see
\cite{Bourgeois:Morse-Bott}, Proposition~2.7):
\begin{equation}\label{index}
\begin{aligned}
\mbox{v--dim}(f):=\mbox{index}(\delbar_f) = &-(2-\overline{p}-\underline{p}
)+\sum_{i=1}^{\overline{p}}
(CZ(\overline{\gamma}_i)+\frac{\dim(\overline{\gamma}_i)}{2})\\
&- \sum_{i=1}^{\underline{p}}
(CZ(\underline{\gamma}_i)-\frac{\dim(\underline{\gamma}_i)}{2}) +
2c_1(T\CC\PP^2)[f].
\end{aligned}
\end{equation}
The  number $c_1(T\CC\PP^2)[u]$ is the Chern number of the bundle
extended into the punctures using the fixed trivialization of it
in the neighborhood $U$ of $L$ given by that of $T(T^\ast L)\to
T^\ast L$. here we use  the symplectic identification of the
neighborhood $U$ of $L$ with a neighborhood of the zero section
$V\supset 0_L$.  The Conley--Zehnder index $CZ$ is given with
respect to the same trivialization.

Notice that all geodesics on the flat surface are minimizing. Thus
their Morse indices $\mbox{index}(\gamma)$ vanish altogether.
Hence the virtual dimension of a holomorphic sphere $f$ with $p$
(negative) punctures in $\CC\PP^2\setminus L$ is given by
\begin{equation}\label{indexsphere}
\mbox{v--dim}(f) = p-2-\sum_{i=1}^p
\mu(\gamma_i)+2c_1(T\CC\PP^2)[f].
\end{equation}

\subsection{Limits of smooth holomorphic spheres}\label{limits}

In this section we describe the behavior of sequences of $J_\tau$--holomorphic
curves as the parameter $\tau$ increases.
We need to describe the objects which will be the limits of
sequences of $J_N$--holomorphic curves in $\CC\PP^2$.

\begin{definition}[Broken holomorphic curves]
A {\em broken $J_\infty$--holomorphic sphere} is a finite
collection of punctured holomorphic curves
$F:=(F^{(1)},...,F^{(N)})$ where

\begin{align}
F^{(1)} &: \dot{S}^{(1)} \longrightarrow X^+\nonumber\\
F^{(k)} &: \dot{S}^{(k)} \longrightarrow \RR\times M \quad \mbox{ for }
k=2,...,N-1\nonumber\\
F^{(N)} &: \dot{S}^{(N)} \longrightarrow X^-\nonumber
\end{align}

and for each {\em level} $S^{(k)}$ denotes the disjoint union of a
finite number of $\CC\PP^1$'s  with finite sets of (pairwise
disjoint) positive and negative punctures
$\overline{z}_1^{(k)},...,\overline{z}_{\overline{p}^{(k)}}^{(k)}$
and
$\underline{z}_1^{(k)},...,\underline{z}_{\underline{p}^{(k)}}^{(k)}$
given on them. As before, we denote by $\dot{S}^{(k)}$ the
corresponding punctured curves.

They  satisfie the following conditions
\begin{enumerate}
   \item $\underline{p}^{(k)} = \overline{p}^{(k+1)}$ with the
         understanding that the first level has no positive and
         the last no negative punctures: $\overline{p}^{(1)} =
         \underline{p}^{(N)} = 0$.
   \item The closed Reeb orbits in $M$ which describe the
         asymptotics  at the punctures agree correspondingly:
         $\underline{\gamma}_j^{(k)} = \overline{\gamma}_j^{(k+1)}$.
   \item Each level is {\em stable}, i.e.~$F^{(k)}$ contains a component
         which is not a cylinder over a closed Reeb orbit without marking.
   \item We compactify
         $\dot{S}^{(k)}$ as described in the definition of punctured holomorphic
         spheres to obtain oriented surfaces with boundaries,
         $\overline{S}^{(k)}$ which are spheres with holes.
   \item A diffeomorphism used for gluing of the boundary
         components of $\overline{\Sigma}^{(k)}$ and $\overline{\Sigma}^{(k+1)}$
         corresponding to $\underline{z}_j^{(k)}$ and $\overline{z}_j^{(k+1)}$,
         is chosen in such a way that it commutes with the restrictions
         of $\overline{f}^{(k)}$ and $\overline{f}^{(k+1)}$ to them.
         There are $m(\gamma_j^{(k)})$ possible choices.
         in the present work we do not care which we choose to formally glue
         the levels.
   \item The surface obtained by gluing all boundary components and denoted by
         $\overline{S}$ is diffeomorphic to the sphere. If we
         assign a graph to $F$ whose nodes correspond to the
         connected components of the $F^{(k)}$ and whose edges to
         the punctures of two adjacent levels which are
         identified, we therefore obtain a tree.
\end{enumerate}
Notice that from the conditions on the asymptotics follows that
$F$ induces a continuous map
$$
\overline{F}:\overline{S}\longrightarrow
\overline{X}:=\overline{X^+}\cup_M \overline{\RR\times
M}\cup_M...\cup_M\overline{X^-}
$$
given by the union of the closures of the $F^{(k)}$. Of course,
the two spaces $\overline{X}\cong X$ are diffeomorphic. Hence the
fundamental class associated to $\overline{F}$, is naturally an
element in the homology of $X$: $[F]\in H_2(X;\ZZ)$.

Finally, notice that the univalent nodes (the branch tips) of the
tree correspond to one--punctured spheres. Since there are no
one--punctured spheres in $T^\ast L$ or $\RR\times U^\ast L$ (see
Lemma~\ref{facts}) these correspond to components of $F^{(1)}$ in
$\CC\PP^2\setminus L$.
\end{definition}

We are now in the position to formulate the splitting theorem for
holomorphic curves in the course of splitting the underlying
almost complex structures $J_{\tau}$  into two almost complex
structures $J^\pm$ as described in Section~ref{neck}. Let ${\cal
J}:=(J^+,J_M,J^-)$ be the triple of almost complex structures.

\begin{proposition}[\cite{BEHWZ, Cieliebak/Mohnke:punc}]\label{compactness}
Let $f_n: \CC\PP^1\to X_{\tau_n}\cong \CC\PP^2$ be a sequence of
$J_{\tau_n}$--holomorphic curves, $\tau_n\to\infty$, with fixed
fundamental class $[f_n]=A\in H_2(X;\ZZ)\cong\ZZ$, equal to the
generator. Then there is a subsequence which converges to a broken
${\cal J}$--holomorphic curve $F=(F^{(1)},...,F^{(N)})$ in the
following sense:

There exists diffeomorphisms $\varphi_n:\overline{S}\to \CC\PP^1$
such that:
\begin{enumerate}
   \item $\lim_{n\to\infty}\varphi_n^\ast j|_{\dot{S}^{(k)}}=j^{(k)}$ in
         $C^\infty_{loc}(\dot{S}^{(k)})$,
         $j^{(k)}$ being
         the conformal structure of the punctured
         holomorphic curve $\dot{S}^{(k)}$ for all $k$.
   \item $\lim_{n\to n}f_n\circ\varphi_n|_{\dot{S}^{(k)}}= F^{(k)}$ in
         $C^\infty_{loc}(\dot{S}^{(k)})$.
\end{enumerate}

Note that it follows that $[F]=A\in H_2(X;\ZZ)$.
\end{proposition}

\begin{remark}
Usually there are two types of possible phenomena for such
sequences of holomorphic curves: non-trivial ''pinching'' of
closed curves (in special situation also known under a different
name, ''bubbling'') as present in Gromov compactness and
''breaking'' apart into different levels of punctured holomorphic
curves similar to the breaking of gradient trajectories in Floer
theory. But if $[f_n]\in H_2(\CC\PP^2;\ZZ)$ are all equal to the
generator there is no pinching. This can be seen as follows.
Assume there is a curve $\gamma$ pinching non-trivially. If we cut
$\CC\PP^1$ open along $\gamma$ we will obtain two holomorphic
curves of are
a greater than some threshold $a_0$. Assuming
$\gamma$ is very small with respect to the metric we may fill
either of the parts by some disk with small symplectic area and
obtain a nontrivial splitting of $A\in H_2(\CC\PP^2;\ZZ)$ into
homology classes which both have positive symplectic area, which
is impossible. Notice that the metric $\omega(.,J_\tau .)$
degenerates as $\tau\to\infty$. In the case that the pinching
approaches the singular locus we have to make use of the fact,
that any holomorphic curve with boundary $\gamma$ either
intersects $\CC\PP^2\setminus U$ or the function $|.|_g$ attains
its maximum on the boundary. If the latter occurs for one of the
two pieces of $\CC\PP^1\setminus \gamma$ it turns out that in the
limit there was no non-trivial pinching. Otherwise we obtain a
contradiction as indicated above.
\end{remark}

We need some more facts about the limit for the proof of the main
result:

\begin{lemma}\label{limitfacts}
\begin{enumerate}
\item With respect to any trivialization of $T(T^\ast L)$ the Chern
      number of $F^{(1)}$ satisfies
      $$
      c_1(F^{(1)}) = f_n^\ast c_1[\CC\PP^1] = 3.
      $$
\item The virtual dimension of $F^{(1)}$ turns out to be
      $$
      \mbox{v--dim}(F^{(1)}) = 4.
      $$
      Equality holds if and only if $N=2$ and $F^{(2)}$
      consists of a disjoint union of $2$--punctured $J^-$--holomorphic
      spheres in $T^\ast L$.
\item The multiplicity of any point $x\in F^{(1)}$ or $x\in F^{(N)}$ is
      $$
      m_x(F) = 1.
      $$
\end{enumerate}
\end{lemma}
\begin{proof}
There is a compact set $K\subset \dot{S}^{(1)}$ such that
$(f_n\circ\varphi_n)^{-1}(\CC\PP^2\setminus U)\subset K$ for all
$n$ sufficiently large. Since $f_n\circ\varphi_n|_K$ converge in
$C^\infty(K)$, the bundle $(f_n\circ\varphi_n|_K)^\ast T\CC\PP^2$
does. On the other hand $(f_n\circ\varphi_n)^\ast T\CC\PP^2$ is
trivialized over $(f_n\circ\varphi_n)^{-1}(U)$ containing all
$\dot{S}^{(k)}$ for $k>1$. By the choice of $K$ we may define
Chern numbers for $(f_n\circ\varphi_n|_K)^\ast T\CC\PP^2$, which
are equal to the Chern number of $f_n$, $f_n^\ast
c_1(T\CC\PP^2)[\CC\PP^1]= 3$ by that remark.

The second statement is a corollary of the first. We find that
\begin{align}
\mbox{v--dim}(F^{(1)}) &= -(2\iota^{(1)} - \underline{p}^{(1)}) -
\sum_{j=1}^{\underline{p}^{(1)}}(CZ(\underline{\gamma}_j^{(1)}) -
\frac{\dim(\underline{\gamma}_j^{(1)})}{2}) + 2c_1(F^{(1)})\nonumber\\
&= 4 - \chi(\overline{F}\setminus F^{(1)}) - \sum_j
\mu(\underline{\gamma}_j^{(1)})\nonumber\\
&=4 - \chi(\overline{F}\setminus F^{(1)}) \le 4\nonumber
\end{align}
Here $\iota^{(1)}$ denotes the number of connected components of
the first level $\dot{S}^{(1)}$.  The second line is due to
Lemma~\ref{viterbo} and the fact, that $\overline{S}\cong S^2$,
the third line uses that the projection to $L$, $\pi\circ
(\overline{F}\setminus F^{(1)})$ provides a $2$--chain in $L$
whose boundary is
$$
\partial(\pi\circ (\overline{F}\setminus F^{(1)})) =
\sum_j [\overline{\gamma}^{(1)}_j],
$$
therefore providing the cancellation of the Maslov indices in the
formula. Finally the last inequality uses the fact that there are
no contractible geodesics in $L$. Therefore the Euler
characteristic of each connected component of
$\overline{F}\setminus F^{(1)}$ is non--positive. Equality holds
therefore if all components of $\overline{F}\setminus F^{(1)}$
consist of annuli.
 But notice that
there the only  holomorphic cylinders in $\RR\times U^\ast L$
(punctured holomorphic spheres with one positive and one negative
puncture), due to the simple fact that closed geodesic (Reeb
orbits) which are homotopic have the same length (action).
Therefore, in this case $\overline{F}\setminus F^{(1)}= F^{(2)}$
consist of $2$--punctured holomorphic spheres in $T^\ast L$.

The last statement follows from the general principle we will use
in the proof of the main result. Fix a non--singular point $x$ on
the top level, and a complex line $\nu\subset T_x\CC\PP^2$. For
any $\tau$ there is a (unique) $J_\tau$--holomorphic sphere
$g_\tau$ through $x$ and tangent to $\nu$ representing the
generator of $H_2(\CC\PP^2;\ZZ)$. For a subsequence of $\tau_n$ of
the compactness statement, we have convergence of the
$g_n:=g_{\tau_n}$ in that sense, to a broken holomorphic sphere
$G$ passing through $x$. Moreover, either the image is tangent to
$\nu$ or it has a cusp singularity at $x$. In both cases the
images of the two punctured holomorphic spheres do not coincide.
Hence there is a neighborhood $W$ of $x$ such that $F\cap G\cap W
= \{x\}$ and they intersect at $x$ with a finite algebraic
intersection index (see \cite{McDuff:intersections,
Micallef/White}), which is equal to the product
$\iota_x(F,G)=m_x(G)m_x(F)$ of the multiplicities of the curves at
$x$. For $n$ large enough the sum of intersection indices of
points in $f_n\cap g_n\cap W$ is equal to that same number. Since
each further point of $f_n\cap g_n$ would contribute a positive
number to the total sum of intersection indices and the algebraic
intersection index of $[f_n]=[g_n]=[H]$ is $f_n\cdot g_n = 1$ we
conclude that
$$
\iota_x(F,G)=m_x(F)m_x(G) = 1.
$$
\end{proof}

We choose the compatible almost complex structure $J^+$ on
$\CC\PP^2\setminus L$ such that all {\em simple} punctured
$J^+$--holomorphic curves are regular. Much less than that is
needed to achieve that (somewhere injective, see
\cite{McDuff/Salamon:Quantum}). In particular, the linearization
of the $\delbar$--equation at each punctured holomorphic sphere in
$\CC\PP^2$ which appears in the limits we consider, is surjective.
That means by the implicit function theorem that the number of
true deformation parameters of it is exactly given by the virtual
dimension. Notice that we do not have to change the structure of
$J^+$ on the ends since there is no non-constant punctured
$J^+$--holomorphic curve in $U\setminus L\subset \CC\PP^2\setminus
L$ by the usual argument using maximum principle (see
\cite{Hofer:Weinstein}).

\section{Proof of Theorem~\ref{threadings}}
We will successively construct two broken holomorphic spheres with
relations to $L$ and to each other which will provide the disks
and spheres we are searching for.

\subsection{Construction of the disks}
Choose a point $x\in L$ which, in the case that $L\cong K$ is the
Klein bottle, lies on one of the odd geodesics. Pick a complex
line $\nu\subset T_x\CC\PP^2$ which is not tangent to $L$:
$\nu\cap T_xL = \{0\}$. Now, for each parameter $\tau$ there is
exactly one $J_\tau$--holomorphic sphere homologous to the
generator of $H_2(\CC\PP^2;\ZZ)$ passing through $x$ and tangent
to $\nu$. Let $f=(f^{(1)},...,f^{(N)})$ denote the broken
holomorphic sphere which is the limit of a subsequence
$f_{\tau_n}$ due to Proposition~\ref{compactness}. Since all
$f_\tau$ pass through $x\in L$ the lowest level $f^{(N)}$ contains
a punctured holomorphic sphere which passes through $x$ and either
has a cusp singularity at $x$ or is tangent to $\nu$. Due to
Lemma~\ref{facts}, (1) and (2) it has to have at least three
(positive) punctures. Pick any one of these punctures and remove
the corresponding edge of the graph corresponding to $f$. Consider
the connected component  of the resulting graph  which does not
contain the node associated to the punctured holomorphic sphere
through $x$. Pick one of the univalent nodes lying in that part.
It corresponds to a one--punctured holomorphic sphere. Due to
Lemma~\ref{facts}, (1), this will be a one-punctured
$J^+$--holomorphic sphere in $\CC\PP^2\setminus L$. We compactify
it to obtain one of the desired disks. Since the punctured
holomorphic curve through $x$ has at least three punctures and due
to Lemma~\ref{limitfacts}, (3), we obtain three disjoint disks
$C,D$ and $E$ with boundary on $L$.

\subsection{Threading the eyes}
We pick a point $y\in C\setminus L$, i.e.~on the image of the
corresponding one--punctured $J^+$--holomorphic sphere. Suppose
that $y$ is a smooth point. Next choose a complex line
$\eta\subset T_y\CC\PP^2$. Assume that there is no punctured
$J^+$--holomorphic sphere in $\CC\PP^2\setminus L$ which appears
as a limit of $J_\tau$--holomorphic spheres homologous to the
primitive class, has less than $4$ deformation parameters and
passes through $y$ and is tangent to $\eta$. Notice that the
tangency makes a priori sense, since $y$ cannot be a singular
point according to the proof of Lemma~\ref{limitfacts}, (3).
Notice also that this gadget can be achieved by choosing
$(y,\eta)$ in general position, since this is a $4$--dimensional
family of parameters to choose from. As in the same proof we
consider the family $g_\tau$ of unique $J_\tau$--holomorphic
spheres passing through $y$ while being tangent to $\eta$. A
subsequence of these converges to a broken holomorphic sphere $g$
with the same property. We claim that it contains only one level
$g=(g^{(1)})$. Therefore, it will be a smooth $J^+$--holomorphic
sphere, i.e.~a mapping
$$
g:\CC\PP^1\longrightarrow \CC\PP^2\setminus L.
$$
Assume on the contrary, that $g=(g^{(1)},...,g^{(M)})$ with $M>1$.
The component $g'$ of $g^{(1)}$ passing through $y$ has  at least
$4$ real deformation parameters. Since all components are simple
and by the choice of $J^+$ thus regular, the virtual dimensions of
each component of $g^{(1)}$ is non--negative. By
Lemma~\ref{limitfacts}, (2), this is only possible if
$$
\mbox{v--dim}(g') = 4,
$$
$M=2$, and $g^{(2)}$ consists of a disjoint union of
two--punctured $J^-$--holomorphic spheres in $T^\ast L$. Notice
that since the index of a component of $g^{(1)}$ is either $0$ or
$4$, it is always even. Hence, due to (\ref{indexsphere}) the Reeb
orbit to which the one--punctured spheres are asymptotic to must
all be odd. It is not difficult to see that therefore all
components of $g^{(2)}$ are of the form $f_\gamma$, for an odd
geodesic $\gamma$. Since all components are simple
(Lemma~\ref{limitfacts}, (3)), the geodesics  have to be simple.
For the same reason each such $f_\gamma$  occurs at most once. On
the other hand, the number of such components of $g$ has to be
even. Indeed, due to Lemma~\ref{facts}, (3), the homological
intersection of $f_\gamma$ with $0_L$ with $\ZZ_2$ coefficients is
$1$. But $L$ could be homotoped away from itself (since there are
nowhere vanishing $1$--forms on $L$). Since $H_2(\CC\PP^2;\ZZ_2)=
\ZZ_2[H]$, $H$ being the projective line, and $H\cdot H = 1$, we
learn that $[L]= 0\in H_2(\CC\PP^2;\ZZ_2)$ (that fact was brought
to our attention by \cite{Nemirovski:kleinbottle}). Therefore,
$[g_\tau]\cdot L = 0$ for all $\tau$ and the above claim follows
since $f$ is the limit of a subsequence of these  $g_\tau$. Hence
both possible, $f_\gamma$ and $f_{\gamma'}$, appear in $g^{(2)}$.
But then $g$ intersects $f$ also in $x$ with a finite positive
algebraic intersection index. Since $g_\tau\cdot f_\tau = 1$ we
obtain a contradiction along the lines of the proof of
Lemma~\ref{limitfacts}, (3). Therefore, our assumption on $M$ was
wrong and we verified the claim. Notice that $g$ may not intersect
neither $D$ nor $E$, again due to the same reasoning using
intersection indices. We call the so-obtained $J^+$--holomorphic
sphere without punctures $H$.

The same procedure applied to $D$ and $E$ yields
$J^+$--holomorphic spheres without punctures $F$ and $G$ with the
desired properties. Notice that since $F,G$ and $H$ are compact
they are  $J_\tau$--holomorphic for sufficiently large $\tau$.

\end{document}